\documentclass[12pt,a4paper]{article}

\newcommand{\expit}{\mathrm{expit}}
\newcommand{\logit}{\mathrm{logit}}
\newcommand{\R}{$\mathsf{R\;}$}
\newcommand{\expect}{\mathsf{E}}
\newcommand{\var}{\mathsf{var}}

\usepackage{color}
\usepackage[affil-it]{authblk}
\usepackage{amsfonts}
\usepackage{amssymb,amsmath}
\usepackage{graphicx}

\begin{document}

\title{Inconsistent treatment estimates from mis-specified logistic regression analyses of randomized trials}

\author{J.N.S. Matthews
  \thanks{\texttt{john.matthews@ncl.ac.uk}: corresponding author}}
\affil{School of Mathematics \& Statistics, Newcastle University,\\
  Newcastle upon Tyne, NE1 7RU, UK\\}

\author{Nuri H. Badi
  \thanks{\texttt{nhb\_2020@yahoo.com}}}
\affil{Statistics Department, University of Benghazi,\\ Benghazi, Libya}  

\date{}
\maketitle

\begin{abstract}
When the difference between treatments in a clinical trial is estimated by a difference in means, then it is well known that randomization ensures unbiassed estimation, even if no account is taken of important baseline covariates.  However, when the treatment effect is assessed by other summaries, e.g. by an odds ratio if the outcome is binary, then bias can arise if some covariates are omitted, regardless of the use of randomization for treatment allocation or the size of the trial.  We present accurate closed-form approximations for this asymptotic bias when important Normally distributed covariates are omitted from a logistic regression.  We compare this approximation with ones in the literature and derive more convenient forms for some of these existing results. The expressions give insight into the form of the bias, which simulations show is usable for distributions other than the Normal.   The key result applies even when there are additional binary covariates in the model.
\end{abstract}

Key words {\it Asymptotic bias; baseline values; logistic regression; probit regression; randomized clinical trial.}

\section{Introduction}

Randomized trials are often analysed using a linear or generalized linear model, so that the treatment effect can be adjusted for important baseline covariates.  However, if some baseline variables cannot be measured, or if their importance is not appreciated, then they will be omitted from the model.  Randomization ensures that the estimate of the treatment effect is unbiassed when relevant covariates are omitted from a linear model.  This is a consequence of the unit-treatment additivity in such models \cite[chapter ~5]{coxplan} and does not necessarily carry over to generalized linear models.  Several non-linear models for which unbiassed estimators are obtained, notwithstanding the omission of covariates, are identified in \cite{gail}, who also show that the important case of binary outcomes analysed using a logistic model is asymptotically biassed when covariates are omitted.

Numerous authors have addressed the problem of the effect of the omission of covariates in logistic regression.  In biostatistical contributions an epidemiological perspective is, perhaps more common \cite{leeecon,hauck,robjew,begglag,dylin,drake}, some authors do focus on randomized trials \cite{gail,neuhaus,neujew,hauckanderson}.  Gail and colleagues \cite{gail} derive approximations for the asymptotic bias in the treatment estimator when all covariates other than the treatment indicator are omitted. The case of two general, scalar, covariates, one of which is fitted and the other omitted is considered in \cite{neujew,drake}.  The main exposition in \cite{neujew} assumes that the covariates are independent but, as the authors explain, this restriction can be relaxed. In all these articles Taylor series approximations are used to provide some indication of the size and direction of the bias, so the expressions derived are necessarily restricted to small parameter values, although whether it is the parameter of the fitted or omitted covariate that needs to be small varies between these contributions. 

In this article we make use of the properties of the extended skew-normal distribution \cite{arnold} and an approximation of the logistic function by the probit to obtain expressions for the least false values \cite[p.25]{claesken} of the fitted covariates when other covariates are omitted.  No use of Taylor series approximations is required, so the expressions give excellent numerical results for a wide range of parameter values and provide useful insight into the form the bias takes in a randomized trial. Our main result applies to a logistic regression with a single binary covariate, which we usually take to indicate the treatment allocation, and an arbitrary number of continuous covariates.  The latter are assumed to follow a multivariate Normal distribution but simulation results show that the results hold for a wider class of covariates.  Explicit forms for the asymptotic bias given in \cite{neuhaus,neujew} are derived for our case and compared with that found using the skew-Normal distribution.  Extensions to allow additional binary covariates is possible, although these extensions require further assumptions. 

In the next section we present the expression for the least false values and in Section \ref{sec:other} related work is explored. Some simulation results are given in Section \ref{sect:sims}, extensions to allow additional binary covariates are discussed in Section \ref{sec:ext} and some conclusions are drawn in the final section.

\section{Least false values}\label{sec:lfvals}

Suppose that the random variable $Y\in \{0,1\}$ is related to a binary covariate $T\in\{-1,1\}$ and further covariates $X_1$ and $X_2$, that have $p$ and $q$ dimensions respectively, by
\begin{equation}\label{eqn:true}
\Pr(Y=1 \mid T,X_1,X_2)=\expit(\mu+\alpha T+\beta_1^TX_1+\beta_2^TX_2)
\end{equation}
where $\expit(u)=\exp(u)/[1+\exp(u)]$.  If the fitted model omits $X_2$, i.e. if
\begin{equation}\label{eqn:fit}
\Pr(Y=1 \mid T,X_1)=\expit(\mu+\alpha T+\beta_1^TX_1)
\end{equation}
is assumed to apply then, as the sample size increases, the maximum likelihood estimates $(\hat\mu, \hat\alpha, \hat\beta_1)$ will tend to the least false values $(\mu^*, \alpha^*, \beta^*_1)$.  From the score equations it can be shown \cite[p25]{claesken} that
\begin{eqnarray}
\expect[\expit(\mu^*+\alpha^*T+\beta^{*T}_1X_1)]&=&\expect[\expit(\mu+\alpha T+\beta_1^TX_1+\beta_2^TX_2)]\label{eqn:expit1}\\
\expect[T\expit(\mu^*+\alpha^*T+\beta^{*T}_1X_1)]&=&\expect[T\expit(\mu+\alpha T+\beta_1^TX_1+\beta_2^TX_2)]\label{eqn:expit2}\\
\expect[X_{1j}\expit(\mu^*+\alpha^*T+\beta^{*T}_1X_1)]&=&\expect[X_{1j}\expit(\mu+\alpha T+\beta^T_1X_1+\beta_2^TX_2)]\label{eqn:expit3},
\end{eqnarray}
where $X_{1j}$ is the $j^\mathrm{th}$ element of $X_1$, $j=1,\ldots,p$ and expectations are taken with respect to the joint distribution of $(T,X_1,X_2)$.

The density of an extended multivariate skew-Normal (ESN) random variable $U\in \mathbb{R}^p$ \cite{arnold} is
\begin{equation}\label{eqn:densSN}
f(u)=\frac{\phi_p(u;\omega,\Omega)\Phi(\zeta^T(u-\omega)+\psi)}{\Phi(\psi/\sqrt{1+\zeta^T\Omega\zeta})},
\end{equation}
where $\zeta$ is a $p$-dimensional parameter, $\psi$ is a scalar, $\phi_p(\cdot; \omega,\Omega)$ is the $p$-dimensional multivariate Normal density with mean $\omega$ and dispersion $\Omega$ and $\Phi(\cdot)$ is the standard Normal distribution function. The mean of the ESN distribution is
\[
\expect(U)=\omega+\frac{\Omega\zeta}{\sqrt{1+\zeta^T\Omega\zeta}}\frac{\phi(\bar\psi)}{\Phi(\bar\psi)},
\]
where $\bar\psi=\psi(1+\zeta^T\Omega\zeta)^{-\tfrac{1}{2}}$ and $\phi(\cdot)=\phi_1(\cdot;0,1)$.

We consider the case when, conditional on $T=t$,  $X=(X_1^T,X_2^T)^T$ follows a multivariate Normal distribution with mean $\nu_t$ and dispersion $\Omega$, $t=-1,1$.  In principle we could allow the dispersion to change with $T$ but analytic progress does not seem possible in this case. We also use $\nu_{t,1}$, $\nu_{t,2}$, $\Omega_{11}$, $\Omega_{22}$, $\Omega_{12}$ and  $\Omega_{21}$ to denote the partition of $\nu_t$ and $\Omega$ induced by the partition of $X$.  If we use the approximation $\expit(u)\approx\Phi(cu)$ with $c=16\sqrt{3}/(15\pi)$ \cite{johnkotz} in (\ref{eqn:expit1}), (\ref{eqn:expit2}) and (\ref{eqn:expit3}), then properties of the ESN distribution provide the approximations
\begin{eqnarray}
\beta^*_1&\approx&\frac{\beta_1+\Omega_{11}^{-1}\Omega_{12}\beta_2}{\sqrt{1+c^2\beta_2^T\tilde\Omega\beta_2}}\label{eqn:alpha1}\\
\mu^*&\approx&\frac{1}{\sqrt{1+c^2\beta_2^T\tilde\Omega\beta_2}}[\mu+\tfrac{1}{2}\beta_2^T\{(\nu_{1,2}+\nu_{-1,2})-\Omega_{21}\Omega_{11}^{-1}(\nu_{1,1}+\nu_{-1,1})\}]\label{eqn:alpha2}\\
\alpha^*&\approx&\frac{1}{\sqrt{1+c^2\beta_2^T\tilde\Omega\beta_2}}[\alpha+\tfrac{1}{2}\beta_2^T\{(\nu_{1,2}-\nu_{-1,2})-\Omega_{21}\Omega_{11}^{-1}(\nu_{1,1}-\nu_{-1,1})\}]\label{eqn:alpha3},
\end{eqnarray}

\noindent where $\tilde\Omega=\Omega_{22}-\Omega_{21}\Omega_{11}^{-1}\Omega_{12}$, is the dispersion of $X_2$ conditional on $X_1$.  Outline details of the derivation can be found in the Appendix.  Note that if $\tilde q=\sqrt{1+c^2\beta_2^T\tilde\Omega\beta_2}$ then the presence of the factor $\tilde q^{-1}$ means that even if $X_1$ and $X_2$ are uncorrelated, then, unlike the linear model, $\beta_1^*\ne \beta_1$ unless, trivially, $\beta_2=0$, or $\tilde\Omega=0$, i.e. the variation in the omitted variables is wholly explained by the fitted variables.  To repeat, the only approximation required for the results in (\ref{eqn:alpha1}), (\ref{eqn:alpha2}) and (\ref{eqn:alpha3}) is that of a logistic by a probit, which is well known to be highly accurate.

When $T$ is the treatment indicator and $X$ are baseline covariates from a randomized trial, then the assumption made above, namely $\var(X\mid T=1)=\var(X\mid T=-1)$ is automatically satisfied and, additionally $\nu_1=\nu_{-1}$, so (\ref{eqn:alpha3}) implies that the least false value of the treatment effect $\alpha$ is
\begin{equation}\label{eqn:mainres}
\alpha^*\approx\frac{\alpha}{\sqrt{1+c^2\beta_2^T\tilde\Omega\beta_2}}.
\end{equation}
Apart from the cases already mentioned which give $\tilde q=1$, (\ref{eqn:mainres}) shows that the omission of relevant covariates means that the treatment estimator will be biassed towards no effect.

\section{Relation with other work}\label{sec:other}

\subsection{No fitted covariates other than the treatment indicator}

Gail and colleagues\cite{gail} considered the bias of treatment estimates for the case when there are no fitted covariates, i.e. the fitted equation is simply
\[
\Pr(Y=1 \mid T)=\expit(\mu+\alpha T),
\]
as opposed to (\ref{eqn:fit}), and where the omitted covariates are not restricted to being Normally distributed.  Finding $\mu^*$ and $\alpha^*$ amounts to solving equations (\ref{eqn:expit1}) and (\ref{eqn:expit2}) with $X_1$ omitted.  In \cite{gail} Taylor series expansions for small $\beta_2^TX_2$ were used to obtain the approximate solution
\begin{equation}\label{eqn:gailbias}
\alpha^*-\alpha \approx -\tfrac{1}{2} \beta_2^T\Omega_{22}\beta_2 (\expit(\mu+\alpha)-\expit(\mu-\alpha)).
\end{equation}

For small $\alpha$ this is approximately $-\alpha\beta_2^T\Omega_{22}\beta_2\expit(\mu)(1-\expit(\mu))$, whereas (\ref{eqn:mainres}) implies that for small $\beta_2^TX_2$ the bias is approximately $-\tfrac{c^2}{2}\alpha\beta_2^T\Omega_{22}\beta_2$.  As $\tfrac{1}{2}c^2\approx 0.173$, this is similar to $\expit(\mu)(1-\expit(\mu))$, which varies from 0.1 to 0.25 as $\mu$ varies over (-2,2).

In \cite{neuhaus,neujew} a different approach was applied to the case when the true model has two scalar covariates, only one of which is included in the fitted model.  As in \cite{gail} no assumption of Normality was made. These authors also used a Taylor series expansion but now applied to the fitted, rather than the omitted covariate.  Using the notation in the present paper, and taking $T$ to be the fitted covariate, the approach in \cite{neujew} noted that 
\begin{equation}\label{eqn:explogit}
\alpha^*=\tfrac{1}{2}[\logit(\pi^*_1)-\logit(\pi^*_{-1})]=H(\alpha)
\end{equation}
where $\pi^*_k=\expect(\expit(\mu+k\alpha+\xi))$, where the expectation is taken with respect to the distribution of $\xi=\beta_2^TX_2$.  Strictly it is the distribution of $\xi$ conditional on $T=k$ but as $T$ is a randomization indicator, this coincides with the unconditional distribution of $\xi$.  Expanding $H(.)$ about $\alpha=0$ \cite{neujew} gives, in the case of logistic regression,  
\begin{equation}\label{eqn:hdash}
\alpha^*\approx \alpha H'(0) =\alpha\; \frac{\pi^*_0-\expect[\expit(\mu+\xi)^2]}{\pi^*_0-\pi^{*2}_0},
\end{equation}
which, as with (\ref{eqn:mainres}), is seen to be closer to 0 than $\alpha$.  Exact analytic evaluation of $H'(0)$ is not possible but further use of the approximation $\expit(u)\approx\Phi(cu)$ and results due to DB Owen reproduced in \cite[p.236]{azzalini}, allow (\ref{eqn:hdash}) to be written as
\begin{equation}\label{eqn:owen1}
\alpha^*\approx \alpha \frac{2T(h,a)}{\Phi(h)\Phi(-h)}=\alpha\frac{T(h,a)}{T(h,1)}
\end{equation}
where $h=c(\mu+\beta_2^T\nu_{,2})/\sqrt{1+c^2\beta_2^T\Omega_{22}\beta_2}$, $a=1/\sqrt{1+2c^2\beta_2^T\Omega_{22}\beta_2}$ and $\nu_{,2}$ is the mean of $X_2$. In (\ref{eqn:owen1}) $T(h,a)$ is Owen's $T$ function \cite{owen}, defined as
\[
T(h,a)=\frac{1}{2\pi}\int_0^a \frac{\exp[-\tfrac{1}{2}h^2(1+x^2)]}{1+x^2}dx,
\]
which has an important role in the computation of bivariate Normal probabilities.  It can be evaluated conveniently by the function \texttt{T.Owen} in the \R package \texttt{sn} \cite{sn}. For fixed $h$, $T(h,a)$ is an increasing function of its second argument and as, in the present application, $0<a<1$, it follows that the expression for $\alpha^*$ in (\ref{eqn:owen1}) is always closer to 0 than $\alpha$. For fixed $a$, $T(h,a)/T(h,1)$ is an even function of $h$ and increases as the magnitude of $h$ increases, so the largest attentuation of $\alpha$ occurs at $h=0$.  As the magnitude of $h$ increases, $T(h,a)/T(h,1)$ approaches one and $\alpha^*$ approaches $\alpha$.

While (\ref{eqn:mainres}) gives a bias in $\alpha$ that does not change with the mean of the covariates, this is not the case with (\ref{eqn:owen1}).  This is most accessibly shown by plotting, for a series of values of $\tilde q^{-1}$, $T(h,a)/T(h,1)$ against $P=\expit(\mu+\beta_2^T\nu_{,2})$, which is a typical reponse probability.  For most randomized trials $P$ will be between 0.1 and 0.9. The figure shows that the bias correction using (\ref{eqn:mainres}) is slightly conservative relative to (\ref{eqn:owen1}) for most values of $P$.  For more extreme $P$, the bias from $T(h,a)/T(h,1)$ reduces, as predicted from the behaviour of this expression for larger $|h|$.

\begin{figure}[!h]
\centerline{\includegraphics[height=4.2in,width=4.2in,angle=0]{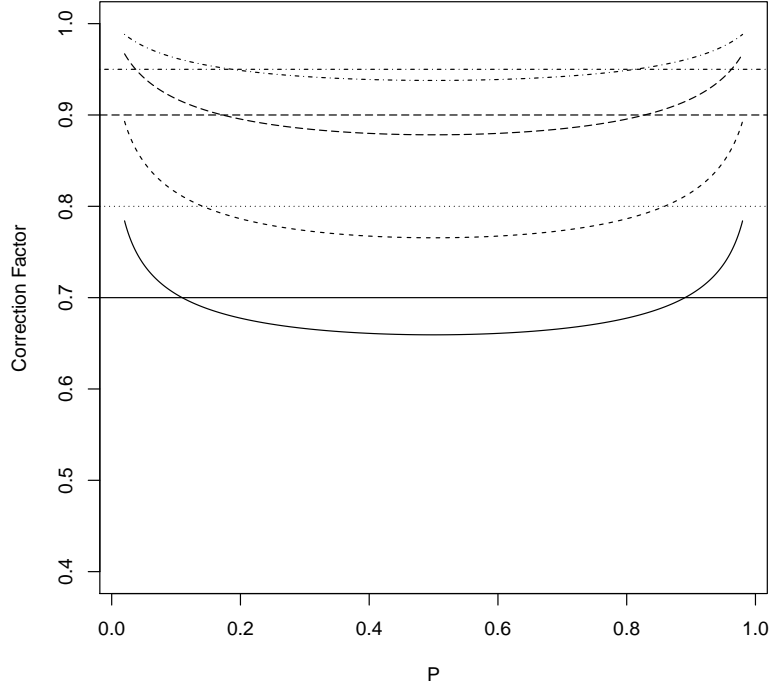}}
\caption{Correction factor $T(h,a)/T(h,1)$ plotted against $P=\expit(\mu+\beta^T_2\nu_{,2})$, for four alternative values of the correction factor $\tilde q^{-1}$, namely 0.7 (solid line); 0.8 (dashed line); 0.9 (long-dashed line); 0.95 (dot-dash line).  The horizontal lines are at the values of $\tilde q^{-1}$}
\label{fig:corrs1}
\end{figure}

\subsection{Covariates fitted in addition to the treatment indicator}

The approach taken in \cite{neujew}, unlike that in \cite{gail}, can be adapted to the case when the fitted model includes covariates $X_1$ in addition to the treatment indicator.  For any given $X_1$ (\ref{eqn:explogit}) still applies, but with the expectation in $\expect[\expit(\mu+k\alpha+\beta_1^TX_1+\xi)]$ now taken with respect to the distribution of $X_2$ given $X_1$.  Consequently the bias factor $T(h,a)/T(h,1)$ still applies but with $a=1/\sqrt{1+2c^2\beta_2^T\tilde \Omega\beta_2}$ and
\begin{equation}\label{eqn:hval}
h=\frac{c[\mu+\beta^T_1X_1+\beta_2^T(\nu_{,2}+\Omega_{21}\Omega_{11}^{-1}(X_1-\nu_{,1}))]}{\sqrt{1+c^2\beta_2^T\tilde \Omega\beta_2}}.  
\end{equation}
This is of limited use because of the dependence on $X_1$, but replacing $X_1$ by its mean $\nu_{,1}$, so  $h=c(\mu+\beta^T\nu)/\sqrt{1+c^2\beta_2^T\tilde \Omega\beta_2}$, provides a workable alternative that can be compared with (\ref{eqn:mainres}) when multiple covariates are fitted.

\subsection{Probit regression}

It is widely acknowledged that in practice logistic and probit regressions can seldom be distinguished in terms of their fit to the data.  As the present analyses have exploited the similarity of $\expit(u)$ and $\Phi(u)$ it is natural to consider the use of probit regression as an alternative to logistic regression, i.e. to replace (\ref{eqn:true}) and (\ref{eqn:fit}) with $\Pr(Y=1 \mid T, X_1, X_2)=\Phi(\mu+\alpha T+\beta_1^TX_1+\beta_2^TX_2)$ etc. The least false values for the maximum likelihood estimators from a probit regression are essentially those in (\ref{eqn:alpha1}), (\ref{eqn:alpha2}) and (\ref{eqn:alpha3}), but with denominator $\sqrt{1+\beta^T_2\tilde \Omega \beta_2}$ in place of $\sqrt{1+c^2\beta^T_2\tilde \Omega \beta_2}$, although the justification of this result is slightly different - see the Appendix for details.  Consequently
\begin{equation}\label{eqn:mainresprobit}
\alpha^* = \frac{\alpha}{\sqrt{1+\beta^T_2\tilde \Omega \beta_2}}
\end{equation}
is an exact expression for the asymptotic bias in the treatment effect from a probit regression with a treatment indicator and Normal covariates. 

Probit regression was also considered in \cite{gail} and \cite{neujew}.  The probit version of (\ref{eqn:gailbias}) is $\alpha^*\approx \alpha(1-\tfrac{1}{2}\beta_2^T\Omega_{22}\beta_2)$.  The bias term given in \cite{neujew} for $H'(0)$ for the probit case is  $\expect(\phi[\Phi^{-1}(\pi_0)])/\phi(\Phi^{-1}[\expect(\pi_0)])$.  If the true model includes both fitted $X_1$ and omitted $X_2$ then the probit analogue of (\ref{eqn:explogit}) applies and the bias factor can be evaluated using  $\pi^*_0=\Phi(\mu+\beta_1^TX_1+\xi)$, with $X_1$ fixed at an arbitrary value and expectations taken over the distribution of $X_2$ conditional on $X_1$.  The denominator of $H'(0)$ is $\phi(h/c)$, with $h$ as in (\ref{eqn:hval}) and the numerator is $\expect[\phi(\mu+\beta_1^TX_1+\xi)]$.  This last expectation has an analytic solution leading to
\[
H'(0)=\frac{\expect(\phi[\Phi^{-1}(\pi_0)])}{\phi(\Phi^{-1}[\expect(\pi_0)])}=\frac{\frac{\phi(h/c)}{\sqrt{1+\var(\xi)}}}{\phi(h/c)}=\frac{1}{\sqrt{1+\beta_2^T\tilde \Omega\beta_2}}.
\]
This coincides with the result from \cite{gail} for small $\beta_2$ and is the same correction factor as obtained from the use of the skew-Normal distribution.  The derivation in \cite{neujew} assumes that $\alpha$ is small and our derivation of the above expression has assumed that the covariates have a multivariate Normal distribution.  In all cases the bias correction for probit regression, unlike logistic regression, depends only on the conditional variance of the omitted variables and their associated regression coefficients, and not on any measure of location.  

\section{Some numerical results}\label{sect:sims}

\subsection{Assessment of the accuracy of the approximations}\label{subsec:normsim}
The simulation results in Table \ref{tab:sims1} assess the accuracy of the forms of $\alpha^*$ for the logistic regression given in equations (\ref{eqn:mainres}), (\ref{eqn:gailbias}) and (\ref{eqn:owen1}), with the last adapted as in (\ref{eqn:hval}) as necessary. The simulated value is found by fitting the reduced model to a sample of size $2\times 10^6$ simulated from the full model: all calculations were performed in \R, version 3.10 \cite{R:2014}.  Three cases are presented: in the first the true model has two Normal covariates, neither of which is fitted, while in the second model only one of these covariates is omitted.  The third model has five covariates, three of which are omitted.  In all cases the Normal covariates have mean 0 and unit variance and correlations are 0.5.  The treatment effect measured by $\alpha$ is taken to be 0.5.  We initially take $\beta_k=0.5$, for $k=1,2$ in the first cases and $1,\ldots,5$ in the final case.  It is important that the simulations correspond to realistic models, with outcome probabilities taking values that are appropriate for a clinical trial.  From (\ref{eqn:true}) we find that 
\[
\Pr(Y=1 \mid T=\pm 1)\approx\Phi\left(\frac{c(\mu\pm\alpha +\beta^T\nu)}{\sqrt{1+c^2\beta^T\Omega\beta}}\right),
\]
so if $\mu$ is chosen so that $\mu +\beta^T\nu=0$ then the outcome probabilities will be around 0.5.

Table \ref{tab:sims1} shows that when $\alpha=0.5$ and $\beta_k=0.5$, all methods perform reasonably when no Normal covariates are fitted, with that from (\ref{eqn:owen1}) doing best.  When some Normal covariates are fitted, Gail's method is not applicable but the proposed extension to (\ref{eqn:owen1}) does well.  The method based on the skew-Normal approximation is conservative, as would be predicted from Figure \ref{fig:corrs1} for response probabilities around 0.5.  When the $\beta_k$ are larger Gail's method fails, as would be anticipated from its derivation. The method due to Neuhaus and colleagues performs better than the skew-Normal factor when no Normal covariates are fitted, but the skew-Normal does better when the fitted model contains some Normal covariates.  The method leading to (\ref{eqn:owen1}) assumes $\alpha$ is small and the final part of Table \ref{tab:sims1} shows that for large $\alpha$ the skew-Normal approximation is again better when Normal covariates are fitted and performs better relative to the method of Neuhaus {\it et al.} than it did for the smaller value of $\alpha$.

\begin{table}[!h]
\begin{center}
\begin{tabular}{l|r|r|r}
\hline
&$p=0,\;q=2$ & $p=1,\;q=1$&$p=2,\;q=3$\\
\hline
&\multicolumn{3}{c}{$\alpha=0.5\;\beta_k=0.5$}\\
\hline
Numerical &0.433&0.482 &0.308  \\
Skew-Normal &0.446& 0.485& 0.330\\
Gail's method&0.408&- & -\\
Neuhaus {\it et al.}&0.434& 0.481& 0.309\\
\hline
&\multicolumn{3}{c}{$\alpha=1.5\;\beta_k=0.5$}\\
\hline
Numerical &1.307&1.447 &  1.328\\
Skew-Normal &1.337&1.454 &1.337 \\
Gail's method&1.262&- & -\\
Neuhaus {\it et al.}&1.302&1.442&1.302 \\
\hline
&\multicolumn{3}{c}{$\alpha=0.5\;\beta_k=2$}\\
\hline
Numerical &0.206&0.347 &0.227  \\
Skew-Normal &0.220& 0.350& 0.220\\
Gail's method&-0.970&- & -\\
Neuhaus {\it et al.}&0.202& 0.330& 0.202\\
\hline
&\multicolumn{3}{c}{$\alpha=1.5\;\beta_k=2$}\\
\hline
Numerical &0.619&1.045 &0.677  \\
Skew-Normal &0.661& 1.051	& 0.661\\
Gail's method&-2.311&- & -\\
Neuhaus {\it et al.}&0.605& 0.990& 0.605
\end{tabular}
\end{center}
\caption{Values of $\alpha^*$ computed using simulation (sample of size $2\times 10^6$) and the three approximations given in equations (\ref{eqn:mainres}), (\ref{eqn:gailbias}) and (\ref{eqn:owen1}), for various values of the regression parameters.  The Normal covariates have mean 0, unit variance and pairwise correlation of $\tfrac{1}{2}$. The number of fitted Normal covariates is $p$ and the number omitted is $q$: throughout $\mu=0$.  }
\label{tab:sims1}
\end{table}

When probit regression is used, the skew-Normal and Neuhaus {\it et al.} expressions coincide and are very close to the simulated value of $\alpha^*$, across a range of values of $\mu$: see Table \ref{tab:sims2}.  Gail's version is reasonable for small $\alpha$ and $\beta_k$ but is poor for larger $\beta_k$.  The lack of dependence of the corrections on $\mu$ is confirmed by the simulated $\alpha^*$, which changes little with $\mu$.  This contrasts with the situation for logistic regression where the simulated $\alpha^*$ show that the bias reduces as $|\mu|$ increases.  The phenomenon applies for all cases but is most clearly seen for small $\beta_k$ and when no Normal covariates are fitted.  This difference between logistic and probit regressions does not appear to be widely appreciated.

\begin{table}[!h]
\begin{center}
\begin{tabular}{l|r|r|r|r|r|r}
\hline
&\multicolumn{6}{c}{$p=0,\;q=2;\;\;\alpha=0.5,\;\beta_k=0.5$}\\
\hline
&\multicolumn{3}{c |}{Logistic regression}&\multicolumn{3}{c}{Probit regression} \\
\hline
& $\mu=0$ &$\mu=2$& $\mu=4$& $\mu=0$ &$\mu=2$& $\mu=4$\\
\hline
Numerical &0.433&0.455 &0.488 &0.378 &0.377& 0.372\\
Skew-Normal &0.446& 0.446& 0.446& 0.378& 0.378& 0.378\\
Gail's method&0.408&0.460 & 0.493& 0.313& 0.313& 0.313\\
Neuhaus {\it et al.}&0.434& 0.452& 0.482& 0.378 & 0.378& 0.378\\
\hline
&\multicolumn{6}{c}{$p=0,\;q=2;\;\;\alpha=0.5,\;\beta_k=2$}\\
\hline
&\multicolumn{3}{c |}{Logistic regression}&\multicolumn{3}{c}{Probit regression} \\
\hline
& $\mu=0$ &$\mu=2$& $\mu=4$& $\mu=0$ &$\mu=2$& $\mu=4$\\
\hline
Numerical &0.207&0.213 &0.231 &0.139 &0.138& 0.140\\
Skew-Normal &0.220& 0.220& 0.220& 0.139& 0.139& 0.139\\
Gail's method&-0.971&-0.139 & 0.390& -2.50& -2.50& -2.50\\
Neuhaus {\it et al.}&0.202& 0.208& 0.227& 0.139 & 0.139& 0.139\\
\hline
&\multicolumn{6}{c}{$p=2,\;q=3;\;\;\alpha=0.5,\;\beta_k=0.5$}\\
\hline
&\multicolumn{3}{c |}{Logistic regression}&\multicolumn{3}{c}{Probit regression} \\
\hline
& $\mu=0$ &$\mu=2$& $\mu=4$& $\mu=0$ &$\mu=2$& $\mu=4$\\
\hline
Numerical &0.437&0.445 &0.458 &0.378 &0.378& 0.379\\
Skew-Normal &0.446& 0.446&0.446 & 0.378& 0.378& 0.378\\
Neuhaus {\it et al.}&0.434&0.452 & 0.482& 0.378 & 0.378& 0.378
\end{tabular}
\end{center}
\caption{Values of $\alpha^*$ computed using simulation (sample of size $2\times 10^6$) and the approximations, for both logistic and probit regression, for different locations of the linear predictor.  The Normal covariates have mean 0, unit variance and pairwise correlation of $\tfrac{1}{2}$. The number of fitted Normal covariates is $p$ and the number omitted is $q$.  }
\label{tab:sims2}
\end{table}

\subsection{Assessment of the effect of departures from Normality}

Some simulations were carried out to assess the effect of non-Normality on the performance of the expressions for $\alpha^*$ in (\ref{eqn:mainres}) and (\ref{eqn:owen1}).  Two types of departure were considered.  The effect of a symmetric non-Normal distribution was assessed by generating $X$ from a central multivariate $t$-distribution with 4 degrees of freedom, while the effect of skewness was assessed using the log-Normal distribution.  In the latter case $X$ was derived from a bivariate Normal variable $W$ with zero mean.  To assess the effect of skewness in the fitted or omitted variable or both, three types of model were considered, with $(X_1,X_2)$ taken as, respectively, $(\exp(W_1)',W_2)$, $(W_1,\exp(W_2)')$ and $(\exp(W_1)',\exp(W_2)')$, where as usual $X_1$ is the fitted covariate and $X_2$ is omitted and $'$ denotes centring to zero mean. The parameters of the $t$ and log-Normal distributions were chosen to give $X_1,X_2$ unit variance and correlation close to $\tfrac{1}{2}$, which implies that the skewness the log-Normal variables are 2.84.  In all simulations $\mu=0$, with $\beta_k=0.5$ or 2 and $\alpha=0.5$ or 1.5, and one scalar covariate is fitted and one omitted.  The correction factors $\tilde q^{-1}$ and $T(h,a)/T(h,1)$ both depend solely on the mean and variance of the $X_i$s, so these will be the same for all of the above models.
\begin{table}[!h]
\begin{center}
\begin{tabular}{l|r|r|r|r}
& \shortstack{$\alpha=0.5$\\ $\beta_k=0.5$} & \shortstack{$\alpha=1.5$\\ $\beta_k=0.5$} & \shortstack{$\alpha=0.5$\\ $\beta_k=2$} & \shortstack{$\alpha=1.5$\\ $\beta_k=2$}\\
\hline
$X$ bivariate $t$, 4 df &0.484&1.456 &0.376 & 1.129\\
$X=(\exp(W_1),W_2)$ &0.479&1.441 &0.352 & 1.061\\
$X=(W_1,\exp(W_2))$ &0.488&1.460 &0.403 & 1.194\\
$X=(\exp(W_1),\exp(W_2))$ &0.481&1.452 &0.375  &1.131\\
\hline
Skew-Normal &0.485& 1.454	& 0.350&1.051\\
Neuhaus {\it et al.}&0.481& 1.442& 0.330&0.990
\end{tabular}
\end{center}
\caption{Values of $\alpha^*$ computed using simulation (sample of size $2\times 10^6$) and the two approximations given in equations (\ref{eqn:mainres}) and (\ref{eqn:owen1}), for various values of the regression parameters.  The covariates have a multivariate $t$ distribution with 4 df or are a mixture of Normal and log-Normal variables.  In each case one covariate is fitted and one omitted, in addition to the treatment indicator:  throughout $\mu=0$. The approximations below the line apply to all the cases above it.}
\label{tab:sims3}
\end{table}

From Table \ref{tab:sims3} we see that for smaller $\beta_k$ the predictions of bias provided by (\ref{eqn:mainres}) and (\ref{eqn:owen1}) remain accurate even when the covariates have non-Normal distributions.  For larger values of $\beta_k$, $\alpha^*$ tends to be closer to $\alpha$ for these non-Normal covariates than for Normal covariates.  However, it should be noted that in this context $\beta_k=2$ is a large coefficient for a covariate with unit variance and unlikely to be encountered in practice.

\subsection{An example: the Mayo Clinic primary biliary cirrhosis trial}

No direct evaluation of the above results is possible as they are all expressed in terms of parameter values.  However, some practical indication of the size of the asymptotic bias, and how this changes with the included covariates, would be helpful.  A trial with binary outcome and several Normal baseline covariates is the primary biliary cirrhosis trial conducted at the Mayo Clinic over ten years from 1974:  the trial randomized patients to placebo or penicillamine \cite{dickson}, and the data are given in \cite{flemharr}. By way of illustration we take end-of-study mortality as the outcome and fit a model with a treatment indicator and five continuous baseline covariates, namely the serum values of bilirubin (mg/dl), cholesterol (mg/dl), albumin (gm/dl), urinary copper ($\mu$g/day) and alkaline phosphatase(AP)  (U/litre).  All variables but albumin were log-transformed (base 10) to achieve Normality.

The dispersion matrix of the five baseline covariates, based on the 312 patients in the trial, was used as $\Omega$ and $\beta$ was taken to be the estimated regression coefficients from the full logistic regression.  The values of $\sqrt{1+c^2\beta_a^T\tilde\Omega\beta_a}$ were then computed for a sequence of models in which the first model includes only the treatment indicator, the second also includes log bilirubin, and then, successively, log cholesterol, albumin and log copper are added.  The correlations are shown in Table \ref{tab:corrs} and the $\tilde q$ values are in Table \ref{tab:taustars}.

\begin{table}
\begin{center}
\begin{tabular}{l|c|ccccc}
&$\beta_k$&log bilirubin& log cholesterol & albumin & log copper&log AP \\
\hline
log bilirubin& 2.19&\emph{0.201} &&&&\\
log cholesterol&-1.43&0.488& \emph{0.036} &&&\\ 
albumin&-0.55 &-0.360 & -0.038& \emph{0.176}&&\\
log copper &0.88&0.598&0.217&-0.278&\emph{0.128}&\\
log AP& 1.69 &0.295&0.351&-0.146&0.277  &\emph{0.098}
\end{tabular}
\caption{The correlations obtained from the dispersion matrix for the five continuous covariates chosen from the PBC trial, with variances on the diagonal and regression coefficients in the second column}\label{tab:corrs}
\end{center}
\end{table}

\begin{table}
\begin{center}
\begin{tabular}{l|c}
Included variables&$\tilde q$\\
\hline
None&  1.311\\
+ log bilirubin &1.072\\ 
+ log cholesterol &1.068\\
+ albumin &1.056\\
+ log copper & 1.039  
\end{tabular}
\caption{The values of $\tilde q=\sqrt{1+c^2\beta_a^T\tilde\Omega\beta_a}$ for a series of increasing models}\label{tab:taustars}
\end{center}
\end{table}
 
If we assume that the model with  treatment indicator and all five variables is the correct model, then $\hat \alpha$ from this model will be asymptotically unbiassed.  However, if a model with no covariates is fitted, $\hat\alpha$ will tend to $\alpha/\tilde q\approx\alpha/1.3$, i.e. a value about 75\% of the correct value.  Including log bilirubin reduces the bias and $\hat\alpha$ will tend to $\alpha/1.07$, a value in error by approximately 7\%.  As Table \ref{tab:taustars} shows, this can be reduced further by including more covariates, although the change is never as marked as when the first variable was introduced.  Of course, different results would be obtained if terms were added in a different order.

\section{Extensions of the model}\label{sec:ext}

The analysis presented thus far applies to a model where, apart from a binary treatment indicator, the covariates are assumed to be continuous. It is often the case that in clinical trials some baseline variables are categorical.  While such variables may have more than two categories, they would usually be included in a linear predictor through dummy variables, so there is no loss in assuming that categorical covariates are binary. The values of the binary treatment indicator are assigned by randomization, so are independent of the values of the other covariates, a feature that would not be shared by a general binary covariate. 

If the model in (\ref{eqn:true}) were extended to include a single non-treatment binary covariate, $B \in \{-1,1\}$, as in
\begin{equation}\label{eqn:c2bincov}
\Pr(Y=1 \mid T,B,X_f,X_a)=\expit(\mu+\alpha T+\gamma B +  \beta_1^T X_1 +\beta_2^T X_2),
\end{equation}
then the foregoing analysis of the effect of omitting $X_2$ from the fitted model can be adapted to this case. In this model, as in Section \ref{sec:lfvals}, $T$ is a binary indicator of the randomized treatment, so is independent of $B$ and $X$. Consequently the parameters defining the distributions of $B$ and $X$ are unaffected by the value of $T$ and we take $\Pr(B=b)=\theta_b$ and $\expect(X \mid B=b)=\nu_b$, $b=-1,1$, but continue to assume that the variance is unaffected by the value of $B$, i.e. $\var(X \mid B=b)=\Omega$.

Under these assumptions it follows that $\beta^*_1$ is as in (\ref{eqn:alpha1}) and
\begin{eqnarray*}
\alpha^*&\approx&\frac{\alpha}{\sqrt{1+c^2\beta_2^T\tilde\Omega\beta_2}}\\
\gamma^*&\approx&\frac{\gamma+\tfrac{1}{2}\beta^T_2([\nu_{1,a}-\nu_{-1,a}]-\Omega_{21}\Omega_{11}^{-1}[\nu_{1,1}-\nu_{-1,1}])}{\sqrt{1+c^2\beta_22^T\tilde\Omega\beta_2}}\\
\mu^*&\approx&\frac{\mu+\tfrac{1}{2}\beta_2^T((\nu_{1,2}+\nu_{-1,2})-\Omega_{21}\Omega_{11}^{-1}(\nu_{1,1}+\nu_{-1,1}))}{\sqrt{1+c^2\beta_2^T\tilde\Omega\beta_2}}.
\end{eqnarray*}
where $\nu_{b,1},\nu_{b,2}$ is the partition of $\nu_b$ corresponding to the partition of $X$ into $X_1$ and$X_2$. The above results are exact for probit regression, provided that the factor $c^2$ is omitted from the denominator.

The above argument can be extended to an arbitrary number of binary covariates, $B_1,\ldots,B_K$ but only at the expense of rather restrictive assumptions about the form of $\expect(X\mid B_1,\ldots,B_K)$.

\section{Discussion}

One of the main reasons for advocating the use of baseline variables in the analysis of a randomized controlled trial is to correct for treatment imbalances in these variables that arise notwithstanding the random allocation.  However, the bias discussed in the present paper is because of the geometrical structure of the statistical model \cite{neujew}, and would apply even if the groups were perfectly balanced.  If there are important covariates in the model then omitting all of them can lead to noticeable biasses, as seen from Table \ref{tab:taustars}, where the log OR is reduced by about 25\%.

Failing to include important covariates can therefore have important consequences for the analysis of a trial.  Not only might the effect of the treatment be underestimated, power calculations may be compromised because observed odds ratios represent less than the true treatment effect.  Of course, in practice the important covariates are not known, at least not with certainty.  If a covariate is included because its importance is suspected, but in fact this view is mistaken, then the corresponding element of $\beta_1$ is zero and there is no penalty in terms of the asymptotic bias in $\alpha$. However, in practical applications, issues related to the finite sample size need to be taken into account. Adding extra variables to a logistic regression may reduce the bias in the treatment estimator but at the expense of an increase in its variance: a fuller investigation of this aspect of the problem is important but beyond the scope of this paper.  Nevertheless, if sufficient data are available when the trial is being planned then it may be possible to use results such as (\ref{eqn:mainres}) to help the triallist make an informed judgment about which covariates ought to be included in the final analysis.

The comparison between logistic and probit analyses is interesting.  If the parameter estimates from a logistic regression are $\hat\beta$ then the estimates obtained from fitting a probit regression to the same data will be approximately $c\hat\beta$, so the corrections in (\ref{eqn:mainres}) and (\ref{eqn:mainresprobit}) are essentially equal.  However, as Figure \ref{fig:corrs1} shows, the correction in (\ref{eqn:owen1}), which can be more accurate than (\ref{eqn:mainres}) for logistic regression, indicates that the asymptotic bias $\alpha^*$ can be greater than is implied by (\ref{eqn:mainres}).  However, the correction in (\ref{eqn:mainresprobit}) is an exact result, so it may be that the problem of asymptotic bias in the estimates of the treatment effect are less if probit is preferred to logistic regression.

\section*{Appendix}

\subsection*{Least false values for logistic regression}

Applying the approximation $\expit(u)\approx \Phi(cu)$ to (\ref{eqn:expit1}), (\ref{eqn:expit2}) and (\ref{eqn:expit3}) and using the properties of the ESN distribution, we obtain from (\ref{eqn:expit1}) and (\ref{eqn:expit2}) the equations
\begin{equation}\label{eqnA:lf12}
p_1\Phi(\psi^*_1)\pm p_{-1}\Phi(\psi^*_{-1})=p_1\Phi(\psi_1)\pm p_{-1}\Phi(\psi_{-1})
\end{equation}
and from (\ref{eqn:expit3}) we obtain
\begin{equation}\label{eqnA:lf3}
\begin{split}
p_1[\nu_1\Phi(\psi^*_1)+\frac{c\Omega_{11}\beta^*_1}{\sqrt{1+c^2\beta^{*T}_1\Omega_{11}\beta^*_1}}\phi(\psi^*_1)]+p_{-1}[\nu_{-1}\Phi(\psi^*_{-1})+\frac{c\Omega_{11}\beta^*_1}{\sqrt{1+c^2\beta^{*T}_1\Omega_{11}\beta^*_1}}\phi(\psi^*_{-1})]\\
=p_1[\nu_1\Phi(\psi_1)+\frac{c(\Omega\beta)_1}{\sqrt{1+c^2\beta^T\Omega\beta}}\phi(\psi_1)]+p_{-1}[\nu_{-1}\Phi(\psi_{-1})+\frac{c(\Omega\beta)_1}{\sqrt{1+c^2\beta^T\Omega\beta}}\phi(\psi_{-1})]
\end{split}
\end{equation}
Here $p_t=\Pr(T=t)$, $\beta^T=(\beta_1^T,\beta_2^T)^T$, $(\Omega\beta)_1$ denotes the first $p$ elements of $\Omega\beta$ and 
\begin{align*}
\psi^*_1&=\frac{c(\mu_1^*+\alpha^*)}{\sqrt{1+c^2\beta_1^{*T}\Omega_{11}\beta^*_1}} &\psi^*_{-1}&=\frac{c(\mu_{-1}^*-\alpha^*)}{\sqrt{1+c^2\beta_1^{*T}\Omega_{11}\beta^*_1}}\\
\psi_1&=\frac{c(\mu_1+\alpha)}{\sqrt{1+c^2\beta^T\Omega\beta}} &\psi_{-1}&=\frac{c(\mu_{-1}-\alpha)}{\sqrt{1+c^2\beta^T\Omega\beta}}
\end{align*}
\noindent
with $\mu^*_t=\mu^*+\beta^{*T}_1\nu_{t,1}$ and $\mu_t=\mu+\beta^T\nu_t$, where $\nu_{t,1}$ is written for the first $p$ elements of $\nu_t$. From (\ref{eqnA:lf12}) we obtain $\psi^*_1=\psi_1$ and $\psi^*_{-1}=\psi_{-1}$, and using this in (\ref{eqnA:lf3}) we get  
\[
\frac{\Omega_{11}\beta^*_1}{\sqrt{1+c^2\beta^{*T}_1\Omega_{11}\beta^*_1}}=\frac{(\Omega\beta)_1}{\sqrt{1+c^2\beta^T\Omega\beta}}
\]
and these can be solved to give (\ref{eqn:alpha1}), (\ref{eqn:alpha2}) and (\ref{eqn:alpha3}).

\subsection*{Least false values for probit regression}

The least false equations for the maximum likelihood estimators for a probit regression differ from (\ref{eqn:expit1}) to (\ref{eqn:expit3}) because of the presence of a weighting factor $\omega=\omega(T,X_1)=\omega(\eta^*)$ with $\eta^*=\mu^*+\alpha^*T+\beta_1^{*T}X_1$, i.e. the $p+2$ equations
\begin{equation}\label{eqnA:probit1}
\expect[\omega(\eta^*)Z\Phi(\eta^*)]=\expect[\omega(\eta^*)Z\Phi(\mu+\alpha T+\beta_1^TX_1+\beta_2^TX_2)]
\end{equation} 
where $\omega(\eta^*)=\phi(\eta^*)/[\Phi(\eta^*)\Phi(-\eta^*)]$, and where $Z$ is taken to be, successively, 1, $T$ and $X_{1j},\;j=1,\ldots,p$.  The presence of $\omega$ means that the skew-Normal distribution cannot be used to evaluate the expectations in the way it was used for logistic regression, but it can be applied to evaluate the right hand expectation in (\ref{eqnA:probit1}) over the distribution of $X_2$ conditional on $T$ and $X_1$, giving
\[
\expect\left[\omega(\eta^*)Z\Phi\left\{\frac{\mu+\beta_2^T(\nu_{T,2}-\Omega_{21}\Omega_{11}^{-1}\nu_{T,1})+\alpha T+(\beta_1+\Omega_{11}^{-1}\Omega_{12}\beta_2)^TX_1}{\sqrt{1+\beta_2^T\tilde\Omega\beta_2}}\right\}\right]
\]
Consequently, if we choose $\beta_1^*$, $\mu^*$ and $\alpha^*$ as in (\ref{eqn:alpha1}), (\ref{eqn:alpha2}) and (\ref{eqn:alpha3}) but with denominator $\sqrt{1+\beta^T_2\tilde \Omega \beta_2}$ as opposed to $\sqrt{1+c^2\beta^T_2\tilde \Omega \beta_2}$, then equations (\ref{eqnA:probit1}) will be satisfied.

\bibliography{rhdraft3pap2}
\bibliographystyle{wileyj}

\end{document}